\documentclass[openany]{book}    
\usepackage{piers}  
\usepackage{subcaption}

\usepackage[
backend=biber,
style=numeric-comp,
maxbibnames=10,
sorting=none,
]{biblatex}

\begin{document}

\newcommand{\bit}[1]{\ensuremath{\textbf{\textit{#1}}}}

\title{Brief research of traditional and AI-based models for IMD2 cancellation}
\maketitle

\author      {A.A. Degtyarev}
\affiliation {MIPT}
\address     {}
\city        {Moscow}
\postalcode  {}
\country     {Russia}
\phone       {}    
\fax         {}    
\email       {degtyarev.aa@phystech.edu}  
\misc        { }  
\nomakeauthor
\author      {N. V. Bakholdin}
\affiliation {MIPT}
\address     {}
\city        {Moscow}
\postalcode  {}
\country     {Russia}
\phone       {}    
\fax         {}    
\email       {bakholdin.nv@phystech.edu}  
\misc        { }  
\nomakeauthor
\author      {A.Y. Maslovskiy}
\affiliation {MIPT}
\address     {}
\city        {Moscow}
\postalcode  {}
\country     {Russia}
\phone       {}    
\fax         {}    
\email       {}  
\misc        { }  
\nomakeauthor
\author      {S.A. Bakhurin}
\affiliation {MIPT}
\address     {}
\city        {Moscow}
\postalcode  {}
\country     {Russia}
\phone       {}    
\fax         {}    
\email       {}  
\misc        {bakhurin.sa@mipt.ru }  
\nomakeauthor

\begin{authors}

{\bf A.A. Degtyarev}, {\bf N. V. Bakholdin}, {\bf A.Y.Maslovskiy} {\bf and S. A. Bakhurin}\\
\medskip
MIPT, Russian Federation\\

\end{authors}

\begin{paper}

\begin{piersabstract}
Due to the limited isolation of duplexer's stopband transceivers operating in frequency division duplex (FDD) encounter a leakage of the transmitted signal onto the receiving path. Leakage signal with the combination of the second-order nonlinearity of the low noise amplifier (LNA) and receiver down-conversion mixer may lead to second-order intermodulation distortion (IMD2) generation thus greatly reducing the receiver sensitivity.

Cancellation of undesirable interferences based on adaptation of traditional models such as memoryless and memory polynomials, spline polynomial based Hammerstein and Wiener-Hammerstein models proved its efficiency in case of well-known nonlinearity nature. 
On the other hand, currently there is an intensive research in the field of nonlinearity detection by means of neural network (NN) structures. NN-based IMD cancellers are effective in the case of unknown interference content due to their high generalization ability. Therefore, NN approach can provide universal model, which is capable of IMD suppression even in case it is hard to separate intermodulation products generated by LNA, down-conversion mixer or even power amplifier in transmitter path. Nevertheless, such structures suffer from high complexity and can`t be implemented in hardware. Current paper presents low-complexity feed-forward NN-based model, which successfully competes with traditional architectures in terms of computational complexity.

The testbench results demonstrate the acceptable performance of provided model, which can be equal to the polynomial nonlinear canceler's performance at a reduced computational cost.

Current paper provides performance and required resources comparison of traditional memory polynomial-based scheme and NN-based model for IMD2 cancellation 

\end{piersabstract}

\psection{Introduction}

During last decades direct-conversion receiver (DCR) has been widely used in modern smartphones because almost all the components of transceiver can be placed within the single chip \cite{dcr_cmos, dcr_saw_less_1, dcr_saw_less_2}. The simplicity of such approach is provided by single-step downconversion. It is very important because modern telecommunication standards require support for multiple frequency bands.

Common approach of self-interference cancellation is introduced by suppression in analog radio frequency (RF) \cite{im2_bluetooth} and/or digital \cite{gebhard_nrls_dct_imd2} domains to bring the interference power to the Rx noise floor. In current paper we focus on the digital second order intermodulation (IMD2) cancellation. Traditional self-interference (SI) cancellation is provided by adaptation of behavioral models, such as memoryless, memory polynomials, spline polynomials etc \cite{behav_model, spline_sic, wiener_hammerst}.


Another promising approach is presented by digital interference cancellation by means of neural network (NN) training \cite{imd2_nn, ffnn, advanced_ml, hlnn}. Modern NN architectures are capable of competing with traditional approaches in terms of performance and computational complexity in digital predistortion (DPD), self-interference cancellation (SIC) tasks for in-band full duplex systems etc \cite{hcrnn, dpd_cnn}. The main property of models proposed for adaptive signal processing tasks is presented by learning of temporal information. For this purpose, NN-based model is fed by the sequence Tx-signal samples. Moreover, input transmitter sequence is divided into several branches with different delays in order to take into account memory effects of the non-linearity.

In current paper we research the IMD2 generated by non-linear distortion of the single RF mixer. Therefore Tx data is introduced by complex data, i.e. quadrature-modulated signal, whereas Rx data is a real-valued sample sequence.

\psection{Traditional and NN methods}\label{SI_cancel_models}


\psubsection{Traditional polynomial models}

Behavioral modelling is introduced by the generalized memory polynomial (GMP) model. Polynomial model is known as a structure, which describes the PA physical properties for different PA kinds and modes~\cite{syrjala2014analysis_oscilator}.

For the task of IMD2 cancellation special case of GMP is decided to be exploited. Moreover, currently we use orthogonal Chebyshev polynomials basis, which could be expressed mathematically as:
\begin{align}
&y_n = \sum_{k=0}^{K-1}\sum_{p=0}^{P-1}\theta_{k,p} T_{p}(\lvert x_{n-d_k}\rvert), 
\label{cheby_1}
\\
&T_{p}(\lvert x_{n-d_k}\rvert)=\text{cos}(n\cdot \text{arccos}(|x_{n-d_k}|)),
\label{cheby_2}
\end{align}

where $\theta_{k,p} \in \mathbb{R}$ - parameters of Chebyshev polynomial model, $x_n$ - samples of baseband (BB) signal, $d_k$ - signal samples delays, $T_p(\cdot)$ - order $p$ Chebyshev polynomial of the first kind.

The coefficients of GMP are searched by least squares (LS) method because model is single layer and linear with respect to parameters. Therefore LS-method guarantess global convergence in a single optimization step \cite{behav_model}.

Nevertheless, in real application cases first-order optimization algorithms are exploited. In this case current model requires optimization parameters fine-tuning, which might be complicated task. Thus, neural network structures are considered below. NN architectures have high generalization ability and, as a result, less sensitive to the optimization parameters tuning.

\psubsection{NN model}

Neural network architecture for IMD2 cancellation is shown in fig.~\ref{fig:nn_scheme}. Since non-linearity is inertial, NN structure is implied to take into account memory effects. 

In current work we realized really small non-linear model based on behavioral modelling approach. As an example architecture to follow we have chosen Wiener-Hammerstein model.
As commonly known Wiener structure contains of cascade memory-less non-linearity and FIR filter functions, in our case we set memory information directly to the non-linearity functions, because the implementation of cascade FIR filters is really expensive according the complexity of it's realization. Thus, memory effects are considered by means of division of input signal into $M$ branches with different delays, as it shown in fig.~\ref{fig:nn_scheme}. Vector of delayed signal samples feeds the sequence of dense layers, which results in the single sample output: 
\begin{align}
    &y_n=\bit{W}_{\text{out}}\sigma_{L-1}(\bit{W}_{L-1}\cdots\sigma_{1}(\bit{W}_{1}\sigma_{0}(\bit{W}_{0}\bit{f}_n))), 
    \label{NN_model_1}
    \\
    &\bit{f}_n=\begin{pmatrix}
        |x_{n-d_0}| & |x_{n-d_1}| & \cdots & |x_{n-d_{M-1}}|
    \end{pmatrix}.
    \label{NN_model_2}
\end{align}
where $M$~-~number of memory delays, necessary to use different streams of signal, $\sigma_{i}(\cdot)$~--~activation function of $i$-th layer, $\bit{W}_{0}\in\mathbb{R}^{K\times M}$, $\bit{W}_{j}\in\mathbb{R}^{K\times K}$, $j\in\overrightarrow{0, L-1}$~--~hidden dense layers weights matrices, $K$ -- number of hidden layer output channels, $L$~--~number of hidden layers, $\bit{W}_{\text{out}}\in\mathbb{R}^{1\times K}$~--~output dense layer, $\bit{f}_n$~--~vector of input delayed signal magnitudes. 
\begin{figure}[!h]
    \centering
        \captionsetup{justification=centering}
    \includegraphics[width=0.85\textwidth]{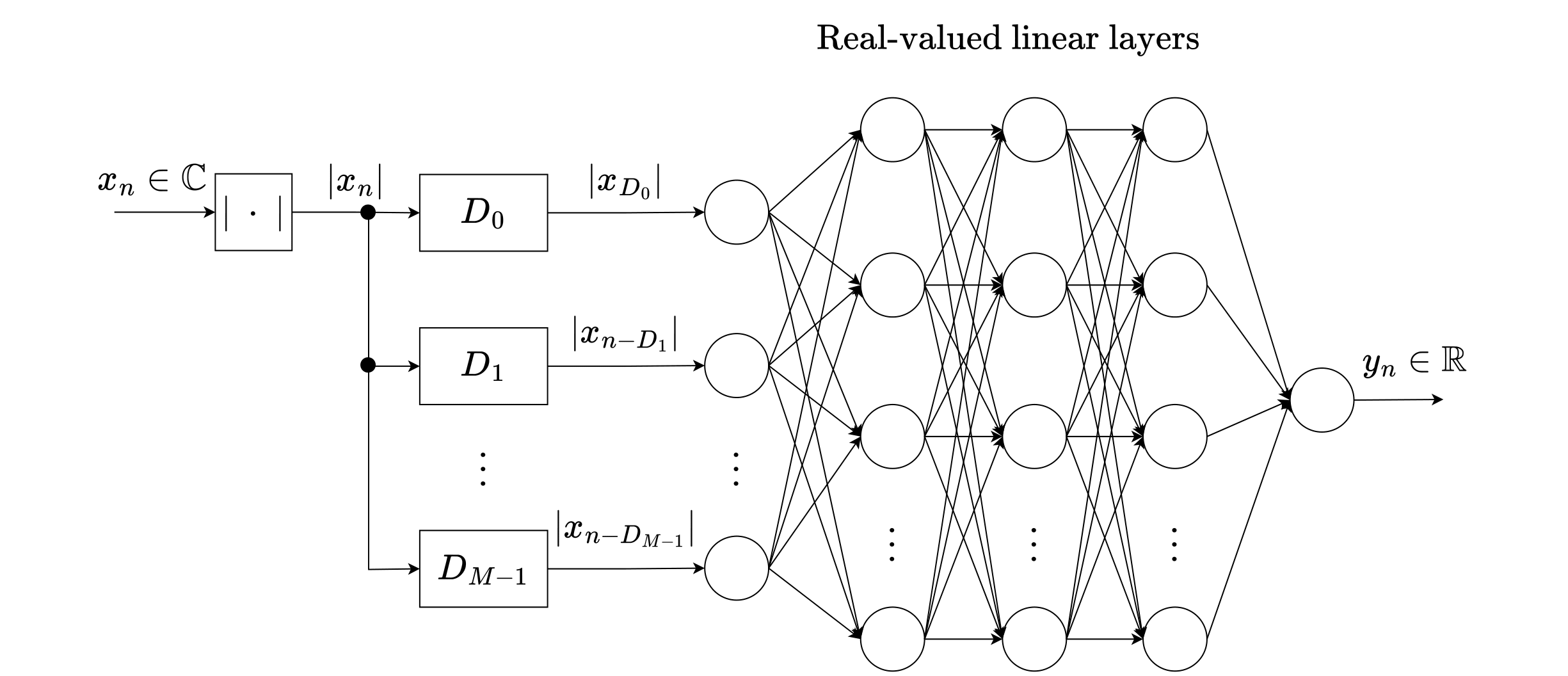}
    \caption{NN-based model architecture}
    \label{fig:nn_scheme}
\end{figure}

Note, that according to the equation~\eqref{NN_model_1} bias was omitted in dense layers for the reason of computational resources conservation.


\psection{Experimental setup}

\begin{figure}[!h]
    \centering
    \captionsetup{justification=centering}
    \includegraphics[width = 0.8\textwidth]{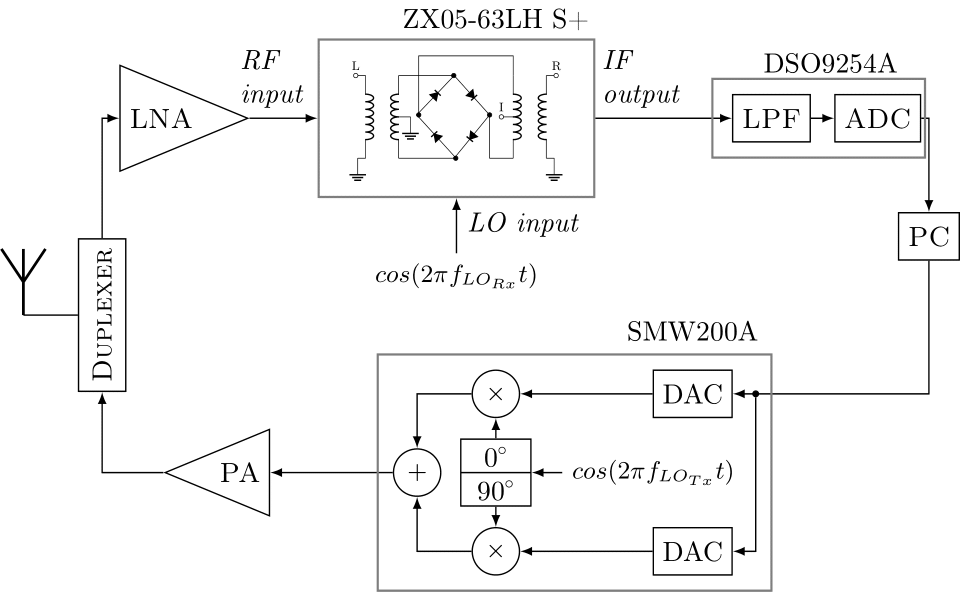}
    \caption{The scheme of testbench}
    \label{fig:testbench_scheme}
\end{figure}

The measurement setup is depicted in fig.~\ref{fig:testbench_scheme} and consists of a personal computer (PC), on which BB data are loaded and then transferred to the signal generator SMW200A. Signal is amplified by nonlinear PA ZRL-3500+, which has a gain 26 dB, 1 dB point $P_{1dB} = 24$ dBm and the OIP3 at 42 dBm. The output of PA is connected with bandpass filter that represents the duplexer with stopband attenuation around 30 dB. The output of BPF is connected to the LNA ZRL-3500+ with NF = 2 dB and gain 26 dB. The ZX05-63LH-S+ level 10 mixer was used for downconversion, it has 37 dB of LO-to-RF isolation and 30 dB of LO-to-IF isolation. No additional analog filter was used in transmit and receive parts.

For the experiments was used complex valued OFDM signal with 5 MHz bandwidth, $f_{Tx}~=~814$ MHz, $f_{Rx} = 859$ MHz, duplex spacing 45 MHz (5G NR Band 26). The baseband (BB) LTE signal is transmitted by R\&S SMW 200A signal generator and amplified by PA. The transmit signal after nonlinear PA leaks through the duplexer stopband into the receiver with 45 MHz frequency offset to the LO signal and is amplified by the LNA. This amplified TxL signal generates the BB IMD2 interference at the output of the mixer. The downconverted signal is captured with digital oscilloscope DSO9254A. The LO signal with the power 10 dBm for the ZX05-63LH-S+ mixer is generated by the R\&S SMW 200A signal generator. 

The transmit power at the output of the PA is set to $P_{Tx}$ = 8 dBm, which in combination with the duplexer attenuation of 30 dB (at $f_{Tx}=814$ MHz) and the LNA gain of 26 dB provides TxL signal power on the input of down conversion mixer $P_{TxL}=8$ dBm $-$ 30 dB $+$ 26 dB $\approx$ 4~dBm.

\psection{Optimization Problem statement}
Based on the described above measurement setup, we need to realize the model, which will be able to compensate leakage and non-linearities of PA problem. To solve this problem there were decided to use the two kind of approaches (section \ref{SI_cancel_models}). Let's denote the to result returned by used model parameterized by the vector $\theta$  at the input x as $ M_{\theta}(x) :=y(x)$. So our problem can be reformulated as problem of supervised learning in the form of regressing task.Then our transmitted signal(TX) and our received signal (RX) be can be used as training dataset(x,$\overline{y}$). Based on the all above we can formulate our problem as minimising of mean squared errors between our RX data and TX data:
\begin{equation} \label{main_problem}
    f(\theta) := \frac{1}{m} \sum_{k=1}^m ([\mathcal{M}_\theta (x)]_k - \overline{y}_k)^2 \rightarrow \min_\theta.
\end{equation}
To assess the quality of the solution obtained as a result of the optimization of this loss functional, we will further use the normalized mean square error quality metric, measured in decibels: 
\[
\label{NMSE}
    \text{NMSE}(y, \overline y) := 10 \log_{10}\left\{\frac{\sum_{k=1}^{m} (y_k - \overline y_k)^2}{\sum_{k=1}^m x_k^2}\right\} \quad \text{dB}.
\]

\psubsection{Least Squares}

It's commonly known that our main problem (\ref{main_problem}) is a least squares problem. For some cases, if model $M_{\theta}(x)$ is linear by coefficients it is possible to find exact solution in this task:

\begin{equation}
    \nabla_{\theta}\left( || A\theta-b||^2+\lambda||\theta||^2)\right)=0
\end{equation}
\begin{equation}
\label{LS solution}
    \theta= (A^H \cdot A +\lambda)^{-1} \cdot A^H \cdot b
\end{equation}
So that, with this way it is possible to find the exect solution for \ref{cheby_1} model, because it is linear by coefficients.
\psubsection{Adam}
According the real simulation , it is impossible to search the solution with as LS. One of the possible way is to use the first order optimization methods. Unfortunately, we can face the problem that Gradient descent method cannot escape from local minima and as a result better performance couldn't be achieved. So that, there were decided to use methods, with 1-st and 2-nd order momentum's, which were proven itself really good for non-convex optimization tasks~\cite{kingma2014adam, maslovskiy2021non}.

\begin{equation}
\label{Adam}
\theta_{k} = \theta_{k-1} - (\frac{\alpha}{\sqrt{\hat v_{k}+ \epsilon}}\hat{m}_{k})\text{,}
\end{equation}
\begin{gather*}
\text{where } m_k = \beta_1 m_{k-1} + (1-\beta_1)\ \nabla f(x) \text{,}\\
v_k = \beta_2 v_{k-1} + (1-\beta_2)\ \nabla f(x)\ ^2
\end{gather*}

\psubsection{L-BFGS}
\label{LBFGS}
As commonly known it is impossible to find least squares solution for non-convex task with respect to model parameters. To achieve the best solution with iterative algorithms we use 2-nd order methods (like Newton).

\begin{equation}
    \theta_k=\theta_{k-1} +\alpha H \nabla_{\theta}f(x)
\end{equation}
\begin{equation*}
    \text{where } H={\left[\nabla^2_{\theta}f\right]^{-1}}
\end{equation*} 

But this algorithm is really heavy, because it is necessary to calculate the inverse Hessian of model. To approximate the inverse second-order derivatives matrix $H$ there were designed a wide range of quasi-Newton approaches. In this regard BFGS algorithm showed itself as robust, high-performant approach. 
\begin{gather}
\label{BFGS}
    \theta_{k+1} = \theta_k - h_k \cdot H_k \nabla f(\theta_k),\text{ where }h_k = \arg \min_{h > 0} f(\theta_k - h \cdot H_k \nabla f(\theta_k)),\\    
    H_{k+1} = H_k + \frac{H_k \gamma_k \delta_k^\top + \delta_k \gamma_k^\top H_k}{\langle H_k \gamma_k, \gamma_k \rangle} - \beta_k \frac{H_k \gamma_k \gamma_k^\top H_k}{\langle H_k \gamma_k, \gamma_k \rangle},\\
    \text{where  }\beta_k = 1 + \frac{\langle \gamma_k, \delta_k \rangle}{\langle H_k \gamma_k, \gamma_k \rangle}, \gamma_k = \nabla f(\theta_{k+1}) - \nabla f(\theta_k), \delta_k = \theta_{k+1} - \theta_k, H_0 = I.
\end{gather}
However, this method is suitable for large-scale problems due to the large amount of memory required to calculate and store the matrix $H_k$ (it is near $n^2$ parameters and iterations). Therefore, in practice the method of recalculating $H_k$ matrix using only $l$ vectors of $\gamma_k$ and $\delta_k$~\cite{zhu1997algorithm, nocedal1980updating}, in this case calculating of $H_{k-l}$ is assumed to be equal to I. the described principle calls L-BFGS($l$) method class with l memory depth.
In our case we used this kind of algorithms to find the best performance really quickly.
\psection{Simulation results}

For NN and polynomial models simulations PyTorch framework was exploited. In current experiments we compared performance and convergence speed of Chebyshev polynomial model~\eqref{cheby_1}, \eqref{cheby_2} where $K=3$, $P=8$ and NN-based canceller~\eqref{NN_model_1}, where $M=3$, $L=2$, $\bit{W}_{0}~\in~\mathbb{R}^{3\times 3}$, $\bit{W}_{1}~\in~\mathbb{R}^{2\times 3}$, $\bit{W}_{\text{out}}~\in~\mathbb{R}^{1\times 2}$. Thus, polynomial model has 24 real parameters, whereas NN-based canceller has 17 real parameters.

Since Chebyshev polynomial is the single-layer model, then its performance can be evaluated by one optimization step of LS algorithm. At the same time, NN structure is a multi-layer model and its possible interference suppression is estimated by L-BFGS algorithm (\ref{LBFGS}).

Commonly in real hardware applications first order algorithms are used for interference cancellation tasks due to the high computational cost of second-order and quasi-Newton methods. Therefore in current simulations convergence speed of NN-based and polynomial-based models trained by gradient descent with Adam optimizer (\ref{Adam}) was compared on the fig.~\ref{fig:lc_400epochs} and fig.~\ref{fig:lc_20000epochs}. 
The power spectrum density (PSD) plots of IMD2 and residual error signals after SIC are shown in the fig.~\ref{fig:psd_plots} and fig.~\ref{fig:psd_plots_magn}.

\begin{table}[h]
    \centering
    \begin{tabular}{||r|r|r|r|r|r|r||}
        \cline{1-7}
        \multirow{2}{*}{Model}&\multirow{2}{*}{algorithm}&\multicolumn{5}{c}{count of iterations} \vline\\ \cline{3-7}
            & & 1000 & 2000 & 5000 & 10000 & 20000 \\
        \hline \hline
      & LS & 23.59 & 23.59 & 23.59 & 23.59 & 23.59 \\
     Polynomial & Adam & 21.96 & 22.40 & 22.76 & 22.83 & 22.91  \\
      & L-BFGS(100) & 23.41 & 23.41 & 23.41 & 23.41 & 23.41 \\
     \hline\hline
      & LS & N/A & N/A & N/A & N/A & N/A \\
     NN & Adam & 21.00 & 21.69 & 22.03 & 23.28 & 23.35  \\
      & L-BFGS(100) & 23.20 & 23.63 & 23.63 & 23.63 & 23.63 \\
        \cline{1-7}
    \end{tabular}
    \centering
    \caption{Comparison of simulation results for polynomial and NN models with different optimizers, suppression vs count of iterations}
    \label{tab:comparison}
\end{table}

\begin{figure}[!h]
    \begin{subfigure}[c]{0.5\textwidth}
        \centering
        \includegraphics[width=\textwidth]{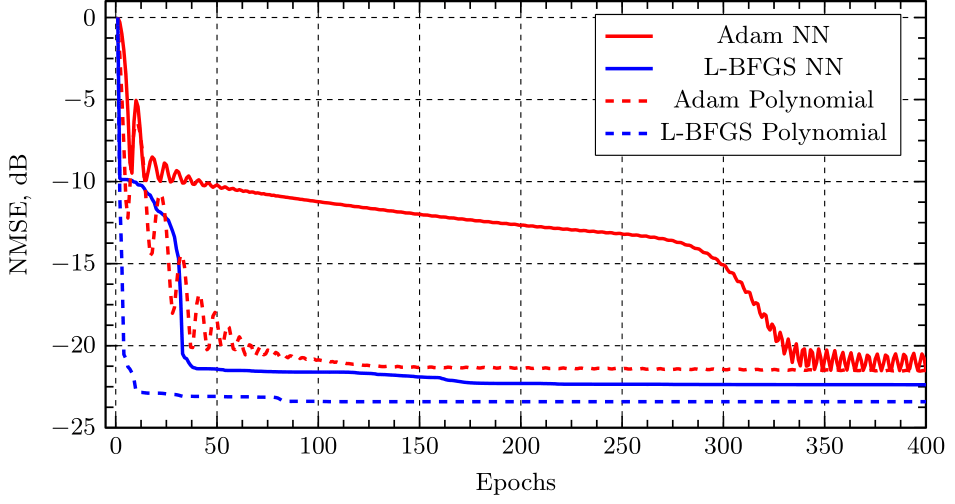}
        \caption{Learning curves for Adam and L-BFGS methods 400 epochs}
        \label{fig:lc_400epochs}
    \end{subfigure}
    \begin{subfigure}[c]{0.48\textwidth}
        \centering
        \includegraphics[width=\textwidth]{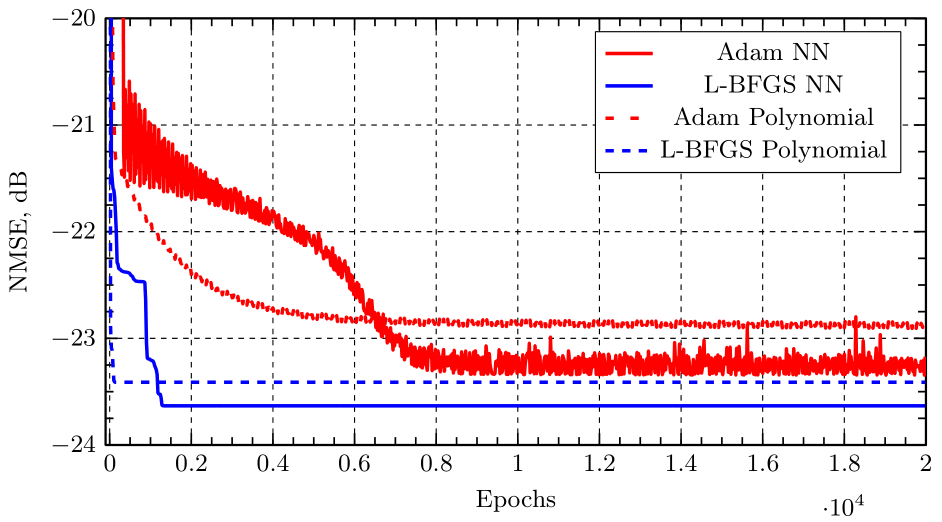}
        \caption{Learning curves for Adam and L-BFGS methods 20000 epochs}
        \label{fig:lc_20000epochs}
    \end{subfigure}
    \caption{Learning curves for Adam and L-BFGS optimizers with different number of epochs}
    \label{learning_curves}
\end{figure}

\begin{figure}[!h]
    \begin{subfigure}[r]{0.5\textwidth}
        \includegraphics[width =\textwidth]{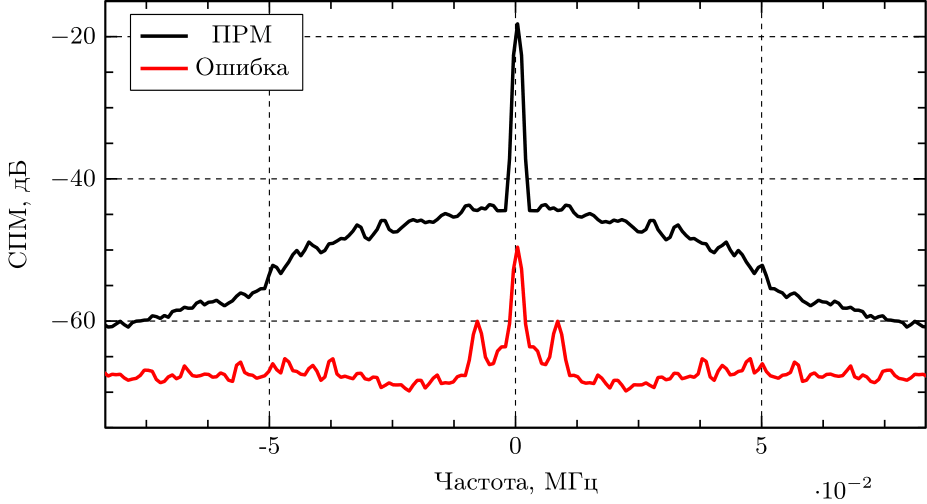}
        \caption{IMD2 before and after cancellation}
        \label{fig:psd_plots}
    \end{subfigure}
    \hfill
    \begin{subfigure}[r]{0.5\textwidth}
        \includegraphics[width =\textwidth]{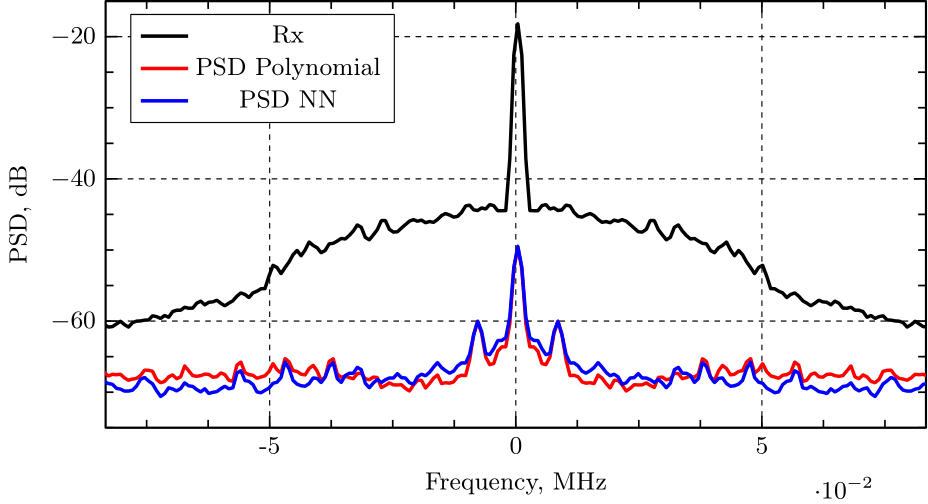}
        \caption{IMD2 before and after cancellation, magnified view}
        \label{fig:psd_plots_magn}
    \end{subfigure}
    \caption{PSD plots of Rx signal, error after cancellation by Chebyshev polynomial and NN models}
    \label{psd}
\end{figure}

According to the results presented in table~\ref{tab:comparison} and fig.~\ref{learning_curves} L-BFGS method provides performance for both architectures close to the LS solution for polynomial NMSE=-23.59 dB. Moreover, for both structures L-BFGS simulation process takes less then 2000 epochs (table~\ref{tab:comparison}).  Current results shows its practical usefulness in terms of models performance evaluation in the field of interference cancellation.

At the same time, first-order method for NN-based model shows higher convergence rate comparing to the polynomial-based canceller due to neural networks generalization ability. For instance, NN architecture achieves 0.44 dB performance improvement in comparison with polynomial in 20000 epochs. Nonetheless, polynomial may achieve full convergence performance by first-order optimizer parameters fine-tuning. This shows one of the remarkable advantages of NN structures.

Although convergence rate of L-BFGS is higher in terms of optimization steps number in comparison with first-order methods, its optimization step time requirements could be higher in contrast to the gradient-based algorithms. Thus, total optimization time of quasi-Newton methods might be comparable to the first-order methods.


\psection{Conclusion}

In the current article we researched NN and polynomial based models for IMD2 cancellation induced by Tx leakage signal in presence of limited stopband attenuation of duplexer.




Current paper presents that both neural network and Chebyshev polynomial based models can achieve good performance but NN model can suppress IMD2 signal without any parameter tuning, whereas for polynomial model requires searching the set of optimal delays.

The findings show that the L-BFGS approach delivers performance for both architectures near to the LS solution for polynomial NMSE=-23.59 dB. Furthermore, the L-BFGS simulation method for both structures requires fewer than 2000 epochs.  Current findings demonstrate its use in the evaluation of models' performance in the interference cancellation domain.

Due to neural networks' capacity for generalization, the first-order technique for NN-based models also demonstrates a greater convergence rate when compared to polynomial-based cancellers. For example, in 20000 epochs, the NN architecture achieves 0.44 dB performance gain over the polynomial. However, polynomial can reach full convergence performance by fine-tuning the first-order optimizer parameters. This demonstrates one of the amazing benefits of NN architectures.

\psection*{Acknowledgment}
The research is supported by the Ministry of Science and Higher Education of the Russian Federation (Goszadaniye), project No. FSMG-2024-0011.

\printbibliography[title = \small REFERENCES]

\end{paper}

\end{document}